\documentclass[reqno]{amsart}
\usepackage{amsmath,amsthm,amssymb,amscd}
\usepackage[pagebackref]{hyperref}

\newcommand{\N}{\mathbb N}
\newcommand{\F}{\mathbb F}
\newcommand{\Z}{\mathbb Z}
\newcommand{\Q}{\mathbb Q}
\newcommand{\R}{\mathbb R}
\newcommand{\C}{\mathbb C}
\newcommand{\eps}{\varepsilon}
\newcommand{\balpha}{\overline{\alpha}}

\newcommand{\bpi}{\overline{\pi}}
\newcommand{\brho}{\overline{\rho}}
\newcommand{\bK}{\overline{K}}
\newcommand{\bQ}{\overline{\Q}}
\newcommand{\cO}{\mathcal O}
\newcommand{\cC}{\mathcal C}

\newcommand{\cT}{\mathcal T}

\newcommand{\fa}{{\mathfrak a}}
\newcommand{\too}{\longmapsto}
\newcommand{\Lra}{\Longrightarrow}
\newcommand{\lra}{\longrightarrow}

\newcommand{\HH}{{\rm H}}
\newcommand{\Qt}{{\mathbb Q}^\times}
\newcommand{\Qts}{\Q^{\times\,2}}
\newcommand{\im}{{\operatorname{im}}\,}
\newcommand{\Gal}{{\operatorname{Gal}}}
\newcommand{\disc}{{\operatorname{disc}}}
\newcommand{\tors}{{\operatorname{tors}}}

\newcommand{\Sel}{{\operatorname{Sel}}}
\newcommand{\Cl}{{\operatorname{Cl}}}
\newcommand{\cl}{{\operatorname{cl}}}
\newcommand{\ts}{\textstyle}
\newfont{\cyr}{wncyb10}
\newcommand{\TS}{\mbox{\cyr Sh}}
\newcommand{\rsm}{\raisebox{0em}[2.3ex][1ex]{\rule{0em}{2ex} }}

\newcounter{lemmacount}[section]

\newtheorem{thm}[lemmacount]{Theorem}
\newtheorem{prop}[lemmacount]{Proposition}
\newtheorem{lem}[lemmacount]{Lemma}
\newtheorem{cor}[lemmacount]{Corollary}

\title[Higher Descent on Pell Conics]{Higher Descent on Pell Conics. \\
      III. The First $2$-Descent}
\author{Franz Lemmermeyer}
\address{Department of Mathematics,
  Bilkent University,
  06800 Bilkent, Ankara, Turkey}
\email{franz@fen.bilkent.edu.tr}

\begin{document}
\maketitle

In \cite{pell1} we have sketched the historical development
of problems related to Legendre's equations $ar^2 - bs^2 = 1$
and the associated Pell equation $x^2 - dy^2 = 1$ with $d = ab$.
In \cite{pell2} we discussed certain ``non-standard'' ideas to 
solve the Pell equation. Now we move from the historical to the 
modern part: below we will describe the theory of the first 
$2$-descent on Pell conics and explain its connections to some 
of the results described in \cite{pell1}, leaving the theory of
the second $2$-descent and its relations to results from \cite{pell2}
to another occasion.

As everyone familiar with the basic arithmetic of elliptic curves
will notice, many of the results (e.g. those on heights) presented
here are special cases of more general theorems. 

\section{Pell Conics}

Since it is our ultimate goal to develop a theory of the
Pell equation that is as close to the theory of elliptic
curves as possible, we will first introduce a more geometric
language. 

We will work over a commutative ring $R$ with a unit element,
which most ofen is $\Z$, $\Z_p$, or a finite field of odd
characteristic. Thus we may and will assume that $R$ is an 
integral domain with a quotient field of characteristic $\ne 2$.

Working with the Pell equation $X^2 - dY^2 = 1$ leads to
numerous problems (not insurmountable, but annoying). For
this reason we will work exclusively with $X^2 - \Delta Y^2 = 4$,
where 
$$ \Delta = \begin{cases}
        d & \text{if}\ d \equiv 1 \bmod 4, \\
       4d & \text{if}\ d \equiv 2, 3 \bmod 4 \end{cases}.$$
Here and in the rest of this article, $d$ will always denote
a squarefree integer; in particular, $\Delta$ is squarefree or
$4$ times a squarefree number. The equation $X^2 - \Delta Y^2 = 4$ 
with $\Delta \in R$ describes a plane algebraic affine curve $\cC$, 
and the set
$$ \cC(R) = \{ (x,y) \in R \times R: x^2 - dy^2 = 4 \} $$
is called the set of $R$-integral points on the conic. 

We now define a group law on the set $\cC(\Q)$ of rational
points on $\cC$ by fixing the neutral element $N = (2,0)$ 
and defining $P + Q = R$ for points $P, Q, R \in \cC(\Z)$ by 
letting $R$ denote the second point of intersection of the 
parallel to $PQ$ through $N$ (see Figure \ref{FAL}).

\begin{figure}[!ht]
[removed]
\caption{Addition Law on Pell Conics}\label{FAL}
\end{figure}

\begin{prop}\label{PGA}
The sum of the two points $P = (r,s)$ and $Q = (t,u)$ in $\cC(\Q)$ is
\begin{equation}\label{Eadd}
 P + Q = \begin{cases}
       \big(\frac{r^2 + \Delta s^2}2, rs \big) 
              = (r^2 - 2, rs) & \text{if}\ P = Q, \\
       \big(2\frac{\Delta (s-u)^2 + (r-t)^2}{\Delta (s-u)^2 - (r-t)^2} \ , \ 
           4\frac{(r-t)(s-u)}{\Delta (s-u)^2 - (r-t)^2} \big) 
                                & \text{if}\ P \ne Q. 
\end{cases} \end{equation}
\end{prop}

Observe that these formulas work in any field in which $\Delta$ 
is not a square; this condition guarantees that the denominator 
$\Delta(s-u)^2 - (r-t)^2$ is nonzero whenever $P \ne Q$.

\begin{proof}
For adding the points $P = (r,s)$ and $Q = (t,u)$, we
have to draw a parallel to the line $PQ$ through $N$
and compute its second point of intersection with $\cC$.
Lines through $N = (2,0)$ have the equation $Y = m(X-1)$.

If $P = Q$, then the slope $m$ of the the tangent at $P$ can 
be computed by taking the derivative of the curve equation 
and solving for $Y'$; we find $Y' = \frac{x}{\Delta y}$, hence 
$m = \frac{r}{\Delta s}$ in $P = (r,s)$. A simple calculation 
yields $X = \frac12(r^2 + \Delta s^2) = r^2 - 2$ and $Y = rs$.

Now assume that $P \ne Q$; if $r = t$, then 
$P = (r,s)$ and $Q = (r,-s)$, and the line through $N$
parallel to $PQ$ is tangent to $N$, that is, we have $P+Q = N$;
this agrees with the formulas above. 

Thus we may assume that $r \ne t$; the line through $PQ$ has slope 
$m = \frac{s-u}{r-t}$. Intersecting this line with $\cC$ leads to
$$ (X-2)\big[X + 2 - \Delta m^2(X-2)\big] = 0; $$
since $X = 2$ gives the point $N$, the $X$-coordinate
of the second point of intersection is given by
$$ X = 2 \frac{\Delta m^2 + 1}{\Delta m^2 - 1}.$$
Plugging in $m = \frac{s-u}{r-t}$, we find
$$ P+Q = \Big(2\frac{\Delta (s-u)^2 + (r-t)^2}
                     {\Delta (s-u)^2 - (r-t)^2} \ , \ 
           \frac{s-u}{r-t} (X - 2) \Big).$$
Now observe that 
$\frac{s-u}{r-t} (X - 2) = 4\frac{(r-t)(s-u)}{\Delta (s-u)^2 - (r-t)^2}$.
\end{proof}

Since we are interested in the integral and not the rational 
solutions of Pell equations, the geometric group law does not
seem to be very helpful. Fortunately, all is not lost:

\begin{prop}\label{PRadd}
The addition formula (\ref{Eadd}) is valid over $\Z$: we have 
$$ 2\frac{\Delta (s-u)^2 + (r-t)^2}{\Delta (s-u)^2 - (r-t)^2} 
           = \frac{rt+\Delta su}2, \quad
   4\frac{(r-t)(s-u)}{\Delta (s-u)^2 - (r-t)^2} = \frac{ru + st}2, $$
hence $P+Q = (\frac{rt+\Delta su}2, \frac{ru + st}2) \in \cC(\Z)$ 
for points $P = (r,s)$ and $Q = (t,u)$ in $\cC(\Z)$.
\end{prop}

\begin{proof}
There is nothing to show if $P = Q$ since, in this case, the coordinates
of $P+Q$ are obviously integral.

Thus we only have to consider the case $P \ne Q$. We have to show 
that the denominator $\Delta (s-u)^2 - (r-t)^2$ divides the numerator. 
Now we can simplify this expression by observing
$$ \Delta (s-u)^2 - (r-t)^2 
   = \Delta s^2 - r^2 + \Delta u^2 - t^2 + 2rt - 2\Delta su 
   = 2(rt - \Delta su - 4). $$
Since $(rt - \Delta su - 4)(ru+st) = 4(r-t)(s-u)$, this gives
$$4 \frac{(r-t)(s-u)}{\Delta(s-u)^2 - (r-t)^2} 
      = 4\frac{(r-t)(s-u)}{2(rt - \Delta su - 4)} 
      = \frac{(rt - \Delta su - 4)(ru+st)}{2(rt - \Delta su - 4)} 
        = \frac{ru + st}2. $$
Observe that if $\Delta \equiv 1 \bmod 4$, then 
$r \equiv s, t \equiv u \bmod 2$, hence $ru+st \equiv 0 \bmod 2$.

Now let us look at the numerator of the $x$-coordinate; since
\begin{align*}
  4(r^2 + \Delta s^2 + t^2 + \Delta u^2) 
   & = (t^2 - \Delta u^2)(r^2 + \Delta s^2) + 
       (r^2 - \Delta s^2)(t^2 + \Delta u^2) \\ 
   & = 2(r^2t^2 - \Delta^2s^2u^2) = 2(rt+\Delta su)(rt-\Delta su),
\end{align*}
we find
\begin{align*}
  2[\Delta (s-u)^2 + (r-t)^2] 
    & = 2[r^2 + \Delta s^2 + t^2 + \Delta u^2  - 2(rt + \Delta su)] \\
    & = (rt+\Delta su)(rt-\Delta su) - 4(rt+\Delta su)] \\
    & = (rt+\Delta su)(rt-\Delta su-4).
\end{align*}
This finally shows
$$ 2 \frac{\Delta(s-u)^2 + (r-t)^2}{\Delta(s-u)^2 - (r-t)^2} 
    = \frac{(rt+\Delta su)(rt-\Delta su-4)}{2(rt - \Delta su - 4)} 
    = \frac{rt + \Delta su}2, $$
and now it follows as before that the x-coordinate of $P+Q$ is integral.
\end{proof}

These addition formulas also show that we have a group law over
any ring in which $2$ is a unit or a prime, such as $\F_q$ for 
odd prime powers $q$, the ring $\Z_p$ of $p$-adic integers and 
its quotient field $\Q_p$, or the rings $\Z_S$ of $S$-integers.

The group law on Pell conics has a well known algebraic 
interpretation: consider the maximal order 
$\cO_K = \Z[\frac12(\Delta + \sqrt{\Delta})\,]$ of the quadratic
number field $K$ with discriminant $\Delta$; sending 
$(x,y) \in \cC(\Z)$ to the unit $\frac12(x+y\sqrt{d}\,) \in \cO_K^\times$ 
induces a bijection $\phi: \cC(\Z) \lra \cO_K^\times$. 

\begin{cor}
The map $\phi$ defined above is an isomorphism of groups. 
\end{cor}

\begin{proof}
Since $\phi$ is bijective, it is sufficient to show that it is 
a homomorphism; but this is clear from 
$$\Big(\frac{r+s\sqrt{\Delta}}2\Big)\Big(\frac{t+u\sqrt{\Delta}}2\Big)
 = \frac12\Big(\frac{rt + \Delta su}2 + \frac{ru+st}2 \sqrt{\Delta}\,\Big)$$ 
and Proposition \ref{PRadd}.
\end{proof}

\section{History of Group Laws}
Describing the history of group laws, whether on elliptic curves
or on conics, is a difficult task for various reasons: first,
because the concept of abstract groups developed very slowly;
in fact, the axioms for abstract groups did not become common 
knowledge until the 1890s. The second reason is that the group 
laws were first discovered in a complex environment: the fact 
that the points on the unit circle $S^1$ form a group had been 
known implicitly since Gauss identified $S^1$ with the set of 
complex numbers with absolute value $1$; these form a group with 
respect to multiplication, as is evident from the relation 
$e^{is} e^{it} = e^{i(s+t)}$ known to Euler. But who first 
realized that the set of {\em rational} points on $S^1$ also 
form a group? 

It is somewhat surprising that the algebraic group structure 
on the unit circle $\cC$ was first defined not over $\Q$ but
over the finite rings $R = \Z/n\Z$: Sch{\"o}nemann \cite{Schm} 
showed that the set $\cC(\Z/n\Z) = \{(x,y) \in \Z/n\Z: 
x^2 + y^2 \equiv 1 \bmod n\}$ is closed with respect to the 
addition $(x,y) + (x',y') = (xx' - yy', xy' + x'y)$.  He also 
showed that $\# \cC(\Z/p\Z) = p - (\frac{-1}p)$ annihilates the 
group $\cC(\Z/p\Z)$. Sch\"onemann's language was algebraic;
the geometric definition of a group law on conics was given by 
Juel \cite[p. 101]{Juel}\footnote{This paper also contains the 
first explicit statement of the group law on elliptic curves.}
who stated it only for circles and hyperbolas. In a review for 
the Fortschritte der Mathematik, St\"ackel \cite{Stae} writes 
about Juel's parametrization of conics (see Figure \ref{FAC}):
\begin{quote}
[Die Parameterdarstellung] beruht auf einer eigent\"umlichen Art
geometrischer Addition, die sich \"ubrigens unter anderem Namen
schon bei v. Staudt findet. Ist n\"amlich $E$ ein fester Curvenpunkt, 
so stehen die drei Curvenpunkte $A, B, C$ in der Beziehung $A+B=C$, 
wenn die Geraden $AB$ und $EC$ sich auf einer festen Geraden $OU$ 
schneiden.\footnote{[The parametrization] is based on a remarkable
way of geometric addition, which can be found in a different guise
already in the work of v. Staudt. In fact, if $E$ is some fixed
point on the curve, then the three points $A, B, C$ on the curve
satisfy $A + B = C$ if and only if the lines $AB$ and $EC$ intersect
on some fixed line $OU$.}
\end{quote}

\begin{figure}[!ht]
[removed]
\caption{Addition Law on Conics}\label{FAC}
\end{figure}

By taking $OU$ to be the line at infinity we recover the geometric
group law defined above.

Veblen \& Young \cite{VY} gave a simplified account of von
Staudt's theory of throws, describing the geometric group law 
on affine lines and on certain conics.

The article \cite{Nie} by Niewenglowski (mentioned by Dickson 
\cite[vol II, p. 396]{Dick}) also contained a hint at the 
geometric group law on conics. Niewenglowski considers the 
hyperbola $x^2 - ay^2 = 1$ and writes
\begin{quote}
Soient $A(1,0)$ le sommet, $A_1(x_1,y_1)$ le premier point 
entier \`a coordonn\'ees positives; la parall\`ele men\'ee 
par $A$ \`a la tangente an $A_1$ donnera le point $A_2(x_2,y_2)$; 
la corde $A_1A_3$ sera parall\`ele \`a la tangente en $A_2$, etc., 
et l'on obtiendra ainsi tous les points \`a coordonn\'ees 
enti\`eres et positives.\footnote{Let $A(1,0)$ be the vertex,
$A_1(x_1,y_1)$ the first integral point with positive 
coordinates; the parallel through $A$ to the tangent at $A_1$
will give the point $A_2(x_2,y_2)$; the secant $A_1A_3$ will
be parallel to the tangent at $A_2$, etc., and in this way we
obtain all the integral points with positive coordinates.}
\end{quote}

The fact that certain arithmetic techniques concerning
curves of genus $1$ admit a geometric interpretation became
common knowledge at the end of the 19th century through the
work of Lucas and Sylvester (see Schappacher \cite{Sch}). 
The algebraic geometer E. Turri\`ere \cite{Tu15} became 
interested in number theoretic problems in 1915, when he 
discussed Fibonacci's question whether $5$ is a congruent 
number using the hyperbolas $y^2-x^2=a$ and $z^2-x^2=b$, 
as well as the cubic $uv(u-v)=av-bu$. In a series of articles
\cite{Tu16,Tu17,Tu18} he then put forward his `arithmogeometry',
a geometric investigation of rational points on algebraic curves.
His plead for a new `arithmetic geometry' seems to have fallen 
on deaf ears; I am not aware of a single reference to these
articles.

Now consider Pythagorean triples $(a,b,c)$, that is, integral
solutions of $a^2 + b^2 = c^2$. We call $(a,b,c)$ primitive
if $\gcd(a,b) = 1$; every Pythagorean triple can be written
in the form $(\lambda a, \lambda b, \lambda c)$ for some
nonzero integer $\lambda$ and a primitive triple $(a,b,c)$, 
and Pythagorean triples that are multiples of the same 
primitive triple are called equivalent.

Idenitfying the equivalence class of the Pythagorean triple 
$(a,b,c)$ with the rational point $(\frac{a}{c}, \frac{b}{c})$ on 
the unit circle gives a group structure to equivalence classes 
of Pythagorean triples. Olga Taussky \cite{Tau} also identified 
the triples $(a,b,c)$, $(-a,b,c)$, $(-b,-a,c)$ and $(a,-b,c)$ 
coming from multiplication by $i$ on $S^1$; thus Taussky's group 
of Pythagorean triples is isomorphic to $\cC(\Q)/\cC(\Q)_\tors$,
where $\cC(\Q)_\tors = \langle (0,1) \rangle$ is the torsion
group of $\cC(\Q)$. Eckert \cite{Eck} proved that this group is 
free abelian, and in fact is a direct sum of infinitely many 
copies of $\Z$, one for each prime $p \equiv 1 \bmod 4$. This 
was rediscovered by Tan \cite{Tan}, who worked with the group 
$\cC(\Q)$ instead. Shastri \cite{Sha} determined the group of 
integral points on the unit circle over number fields.

Other articles dealing with group (or ring) structures on the
set of Pythagorean triples are Baldisserri \cite{Bal},
Beauregard \& Suryanarayan \cite{BS1,BS2,BS3}, Dawson \cite{Daw}, 
Grytczuk \cite{Gry}, Hlawka \cite{Hla}, Wojtowicz \cite{Woi}, and 
Zanardo \& Zannier \cite{ZZ}, whereas Mariani \cite{Mar} and Morita 
\cite{Mor} study groups acting on Pythagorean triples.

In the modern mathematical literature, the group law on conics 
is hardly ever discussed; an exception is the book \cite{PrS}
by Prasolov \& Solovyev, or the web site 
\begin{center}
 {\tt http://www-cabri.imag.fr/abracadabri/Algebre/Groupes/FoliumD.html},
\end{center}
which contains a detailed exposition of the group law on conics.

\section{The First $2$-Descent}

The conic $\cC: X^2 - \Delta Y^2 = 4$ comes attached with an isomorphism
$$\psi: \cC(\Q) \lra K^\times[N]: (x,y) \too \frac{x+y\sqrt{\Delta}}2 $$ 
from the group of rational points on $\cC$ to the elements of norm
$1$ in $K^\times$, where $K = \Q(\sqrt{\Delta}\,)$ is the quadratic
number field with discriminant $\Delta$. We know that $\psi$ restricts 
to an isomorphism $\cC(\Z) \lra \cO_K^\times$. 

\subsection{The Set of First Descendants}\label{SD1}
Now consider any integral point $(x,y) \in \cC(\Q)$ on the Pell conic 
$\cC: X^2 - \Delta Y^2 = 4$. Write $\Delta y^2 = x^2 - 4 = (x-2)(x+2)$. 
Since $\gcd(x+2,x-2) \mid 4$, there are three possible cases:
\begin{enumerate}
\item $x \equiv 1 \bmod 2$: then $\Delta \equiv 5 \bmod 8$,
      $\gcd(x-2,x+2) = 1$, hence $x+2 = ar^2$ and $x-2 = bs^2$, 
      where $ab = \Delta$. Thus $ar^2 - bs^2 = 4$. 
\item $x \equiv 2 \bmod 4$: then we find $\gcd(x-2,x+2) = 4$, 
      hence $x+2 = ar^2$, $x-2 = bs^2$, and again $ar^2 - bs^2 = 4$.
\item $x \equiv 0 \bmod 4$: then $\Delta = 4d$ with $d \equiv 3 \bmod 4$
      and $\gcd(x-2,x+2) = 2$, so $x+2 = 2Ar^2$, $x-2 = 2Bs^2$ with 
      $ab = d$, hence $ar^2 - bs^2 = 4$ for $a = 2A$, $b = 2B$ and 
      $ab = \Delta$.
\end{enumerate}

The curves $\cT_a:ar^2 - bs^2 = 4$ are called the first 
descendants of $X^2 - \Delta Y^2 = 4$. Every integral point
on $\cC$ comes from an integral point on one of the descendants.

If $\Delta < 0$, then $x^2 + |\Delta| y^2 = 4$ implies 
that $x^2 \le 4$, which in turn shows that $x+2 > 0$ 
unless $x = -2$. Thus the descendants all have the form 
$\cT_a$ for positive integers $a$.

If $\Delta > 0$, then $x \ge 2$ or $x \le -2$. The points with 
$x > 0$ come from descendants $\cT_a: ar^2 - bs^2 = 4$ with $a > 0$. 
If $(x,y)$ is such a point, then the points $(-x, \pm y)$ will come 
from the descendant $\cT_{-\Delta/a}: -br^2 + as^2 = 4$ (or, if
$4 \mid \Delta$, from $\cT_{-\Delta/4b}$) describing the same curve
(up to a change of variables) as $\cT_a$. It is therefore 
sufficient to consider descendants $\cT_a$ for $a > 0$ squarefree. 

\begin{thm}
Every integral solution $(x,y)$ of the Pell equation 
$X^2 - \Delta Y^2 = 4$ gives rise to an integral solution of one 
of the equations $\cT_a: ar^2 - bs^2 = 4$, where $a$ and $b$ are 
integers such that $ab = \Delta$, and where $a$ is squarefree. 

Conversely, any integral solution $(r,s)$ of $\cT_a$ gives rise
to an integral solution $(x,y)$ of the Pell equation, where 
$x = ar^2-2$ and $y = rs$. 
\end{thm}

\medskip\noindent{\bf Remark 1.}
If $\Delta = 4d$ with $d \equiv 3 \bmod 4$, the descendants
$\cT_{2a}: 2ar^2 - 2bs^2 = 4$ with $4ab = \Delta$ coincide
with the curves $ar^2 - bs^2 = 2$ occurring in the theory of
Legendre (see \cite[Section 2]{pell1}).

\medskip\noindent{\bf Remark 2.}
Assume that $ar^2 - bs^2 = 4$, where $ab = \Delta$. If 
$s = 1$ is a solution, then $b = ar^2 - 4$, hence 
$\Delta = ab = a(ar^2-4) = a^2r^2 - 4a$. A solution
$s = 2$ implies that $r = 2m$ and leads to $\Delta = a^2m^2 -a$. 
Similarly, solutions $r = 1, 2$ leads to values of $\Delta$ that
are of Richaud-Degert type $\Delta = n^2 + r$ with $r \mid 4n$.

\medskip\noindent {\bf Example.}
Consider $\cC(\Z)$ for $\cC: x^2 - 205y^2 = 4$. The associated 
descendants with an integral point $(r,s)$ and the corresponding 
point $(x,y)$ on $\cC$ are given below:
$$ \begin{array}{r|r|c|c}
  a & \cT_a(\cC) \quad  & (r,s) & (x,y) \\ \hline
\rsm   1 &    r^2 - 205s^2 = 4 & (2,0) & (2,0) \\
\rsm   5 &   5r^2 -  41s^2 = 4 & (3,1) & (43,3) \\
\rsm  41 &  41r^2 -   5s^2 = 4 & (\frac13,\frac13) 
                 & (\frac{23}9,\frac19) \\
\rsm 205 & 205r^2 -    s^2 = 4 & (\frac23,\frac{28}3) 
                 & (\frac{802}9,\frac{56}9)
   \end{array} $$
The existence of integral points on the last two 
descendants cannot be excluded via congruences alone;
this is a case where a second $2$-descent would help.

\subsection{The Group Structure}

The number of descendants we have to consider is always a power of 
$2$, as is the number of descendants with an integral point. This 
could be explained by giving this set of descendants the structure 
of an elementary abelian $2$-group. How can we accomplish this?
 
\subsubsection*{1. The Naive Construction}
The naive idea is to make the first descendants into an
elementary abelian group by defining $\cT_a \cdot \cT_b = \cT_c$,
where $ab = cm^2$ for integers $c, m$ with $c$ squarefree. 
This is easily seen to coincide with the group structure
defined by Dickson \cite[\S 25]{DiSN} (see also \cite{pell1}).

\subsubsection*{2. Using the Group Structure on the Pell conic}
The set of descendants with a rational point can be given
a group structure as follows:  given $(r,s) \in \cT_a(\Q)$ 
and $(t,u) \in \cT_b(\Q)$, compute the corresponding rational
points $(x,y)$ and $(z,w)$ on the Pell conic; the sum 
$(x,y) + (v,w)$ on $\cC(\Q)$ will then come from a rational 
point on some descendant $\cT_c$, and we put
$\cT_a \oplus \cT_b = \cT_c$.

\bigskip

In order to decide whether these group laws coincide (on the
subset of descendants with a rational point) or not, we need 
a better way of finding the descendant $\cT_a$ to which an 
$(x,y) \in \cC(\Q)$ gives rise. Observe that since 
$x = ar^2-2$, we can recover $a$ by mapping $(x,y) \in \cC(\Z)$ 
to the coset $(x+2)\Qts = a \Qts$. Actually, we get a mapping 
$\alpha: \cC(\Q) \lra \Qt/\Qts$ by putting 
$\alpha(x,y) = (x+2)\Qts$ for all $(x,y) \ne (-2,0)$;
using the equation $x^2-4 = \Delta y^2$, we see that we have 
$(x+2)\Qts = (x-2)\Delta \Qts$ whenever both sides are defined, 
and this suggests we define $\alpha(-2,0) = -\Delta \Qts$.

\begin{prop}
Define a map $\alpha: \cC(\Q) \lra \Qt/\Qts$ by
$$ \alpha(x,y) = \begin{cases}
                (x+2)\Qts     & \text{if}\ x \ne -2, \\
                -\Delta \Qts  & \text{if}\ x = -2. \end{cases} $$
If $P = (x,y) \in \cC(\Z)$ with $x > 0$, then $P$ gives rise 
to an integral point on the descendant $\cT_a(\cC)$, where $a$
is a positive squarefree integer determined by 
$\alpha(P) = a\Qts$.                 
\end{prop}

\subsection{The Weil Homomorphism}
The map $\alpha: \cC(\Q) \lra \Qt/\Qts$ is a map between two
abelian groups; is it a homomorphism? Before we show that the
answer is yes, we will give another way to motivate the 
definition of $\alpha$. 

Consider the Pell conic $X^2 - \Delta Y^2 = 4$. We want to define 
a `Weil homomorphism' $\alpha: \cC(\Q) \lra \Qt/\Qts$ with 
kernel $\ker \alpha = 2\cC(\Q)$. Since $2(r,s) = (r^2 - 2,rs)$, 
we could try to map $(x,y)$ to the coset $(x+2)\Qts$; this defines 
a map annihilating $2\cC(\Q)$, but is not defined for $P = (-2,0)$.
On the other hand, we also have $2(x,y) = (2 + \Delta y^2,xy)$; the 
map $(x,y) \too \Delta(x-2)\Qts$ is defined
except for $(x,y) = (2,0)$, and it agrees with the map
defined before for all points $\ne (\pm 2, 0)$. 

Now we claim

\begin{thm}\label{Thom}
The map $\alpha: \cC(\Q) \lra \Qt/\Qts$ is a group homomorphism.
\end{thm}

This will be proved using Galois cohomology below. Before we do 
this, let us derive a few consequences.

\begin{cor}
The group laws defined on the set of first descendants coincide.
\end{cor}

\begin{proof}
Assume that the points $P$ and $Q$ on $\cC(\Q)$ give rise to points 
on the descendants $\cT_a$ and $\cT_b$; then $\alpha(P) = a\Qts$,
$\alpha(Q) = b\Qts$, and since $\alpha$ is a group homomorphism,
$\alpha(P+Q) = ab\Qts$, hence $P+Q$ gives rise to a point on the
descendant $\cT_c$ with $ab = cm^2$ and $c$ squarefree.
\end{proof}

\begin{prop}
The image of $\alpha: \cC(\Z) \lra \Qt/\Qts$ consists of all square 
classes $a\Qts$ for which $ab = \Delta$ for $a, b \in \Z$ and 
$ar^2 - bs^2 = 4$ has an integral solution.
\end{prop}

\begin{proof}
If $a \Qts \in \im \alpha$, then there is a $P = (x,y) \in \cC(\Z)$
such that $\alpha(P) = a\Qts$, and by our construction above the
point $P$ comes from an integral point on $ar^2 - bs^2 = 4$. 
Conversely, if $ar^2 - bs^2 = 4$ has an integral solution, then it 
gives rise to the integral point $P = (ar^2-2,rs)$ on the associated
Pell conic, and $\alpha(P) = (x+2)\Qts = a\Qts$.
\end{proof}

This shows

\begin{cor}
The image of $\alpha: \cC(\Z) \lra \Qt/\Qts$ is finite.
\end{cor}

\begin{proof}
This follows at once from the observation that there are only
finitely many classes $a\Qts$ with $ab = \Delta$ and $a, b \in \Z$.
\end{proof}

Now we claim

\begin{thm}
We have an exact sequence 
$$ \begin{CD} 
   0 @>>> 2\cC(\Z) @>>> \cC(\Z) @>{\alpha}>> \Qt/\Qts. 
   \end{CD} $$
\end{thm}

\begin{proof}
We claim that the kernel of the homomorphism 
$\alpha: \cC(\Q) \lra \Qt/\Qts$ is $\ker \alpha = 2\cC(\Q)$. 
Moreover, the kernel of the induced map
$\cC(\Z) \lra \Qt/\Qts$ is $2\cC(\Z)$.

One direction is clear: if $(x,y) = 2(r,s)$ for some $(r,s) \in \cC(\Q)$,
then $x = r^2 - 2$, hence $x+2 = r^2$ is a square, and this means that 
$(x,y) \in \ker \alpha$.

For the converse, observe that $(x,y) \in \ker \alpha$ if and only if 
$x+2 = r^2$ for some $r \in \Q$. 
Next, $\Delta y^2 = x^2-4 = (x-2)(x+2)$, hence 
$\Delta y^2 = (x-2)r^2$, and thus $x-2 = \Delta s^2$ for some 
$s \in \Q$. On the other hand, $x-2 = x+2-4 = r^2 - 4$, 
hence $r^2 - \Delta s^2 = 4$. Thus $(r,s) \in \cC(\Q)$, and 
it is easily checked that $2(r,s) = (x,y)$.

Now consider the restriction of $\alpha$ to $\cC(\Z)$. If $x \in \Z$
in the above proof, then clearly $r \in \Z$, and $r^2 - \Delta s^2 = 4$
then implies that we also have $s \in \Z$ if $(x,y) \in \cC(\Z)$.
\end{proof}

This immediately implies

\begin{cor}[Weak Theorem of Mordell-Weil]
The group $\cC(\Z)/2\cC(\Z)$ is finite.
\end{cor}

In the next section, we will use the theory of heights to prove
that $\cC(\Z)$ is finitely generated. This implies that  
$\cC(\Z) \simeq \cC(\Z)_\tors \oplus \Z^r$ for some $r \ge 0$, 
and the fact that the torsion group $\cC(\Q)_\tors$ is cyclic 
shows that $\cC(\Z)/2\cC(\Z) \simeq (\Z/2\Z)^{r+1}$. 
This is the analog of Tate's formula for the
$2$-rank of an elliptic curve with rational $2$-torsion:

\begin{prop}
We have $\cC(\Z) \simeq \cC(\Z)_\tors \oplus \Z^r$, where $r \ge 0$
is determined by $\im \alpha = 2^{r+1}$.
\end{prop}

This also implies

\begin{thm}\label{TDi}
Consider the Weil map $\alpha> \cC(\Z) \lra \Qt/\Qts$ for the Pell 
conic $\cC: X^2 - \Delta Y^2 = 4$, where $\Delta > 0$. The following 
assertions are equivalent:
\begin{enumerate}
\item $\cC(\Z) \simeq \Z/2\Z \oplus \Z$;
\item $\#\ \im \alpha = 4$.
\end{enumerate}
\end{thm}

The implication $(1) \Lra (2)$ of Theorem \ref{TDi} is a modern
formulation of Dirichlet's Theorem \cite[Thm. 3.3.]{pell1}.

\subsection*{Proof of Theorem \ref{Thom}.}
Let $\cC: X^2 - dY^2 = 4$ denote the Pell conic, and 
$[2]: \cC(K) \lra \cC(K)$ multiplication by $2$. 

\begin{prop}\label{Pdiv2}
We have an exact sequence
\begin{equation}\label{ES0} \begin{CD}
   0 @>>> \cC(\bQ)[2] @>>> \cC(\bQ) @>[2]>> \cC(\bQ) @>>> 0,
   \end{CD} \end{equation}
where $\cC(\bQ)[2] = \{(-2,0), (2,0)\} = \cC(\Q)[2]$.
\end{prop}

\begin{proof}
Let us first prove that $[2]$ is surjective. Given 
$(r,s) \in \cC(\bQ)$, we find that $2(x,y) = (r,s)$ implies 
$r = x^2 - 2$ and $s = xy$. Thus $x^2 = r+2$, and either $y = 0$ 
(if $r = -2$) or $y = \frac{s}{x}$. In either case, 
$(x,y) \in \cC(\bQ)$ satisfies $2(x,y) = (r,s)$. 

The same formulas show that $\ker [2] = \{(\pm 2, 0)\}$: in fact,
if $(r,s) = (2,0)$, then $x^2 = r+2 = 4$ implies $x = \pm 2$
and $y = 0$.
\end{proof}

Now let $G = \Gal(\bQ/\Q)$ denote the absolute Galois group of $\Q$.
Since $\cC(\bQ)[2]$ consists of rational points, we have
$\cC(\bQ)[2] \simeq \Z/2\Z$ as Galois modules, and the long exact 
cohomology sequence gives
\begin{equation}\label{ES1} \begin{CD}
    \cC(\Q) @>[2]>> \cC(\Q) @>>> \HH^1(\Z/2\Z) @>>> \HH^1(\cC) @>[2]>> \HH^1(\cC),
   \end{CD} \end{equation}
where $\HH^1(A) = \HH^1(G,A)$ and $\cC = \cC(\bK)$. 

Next we compute $\HH^1(\Z/2\Z)$; we start with the Kummer sequence 
$$ \begin{CD}   
     1 @>>> \Z/2\Z @>>>  \bQ^\times @>[2]>> \bQ^\times @>>> 1 
   \end{CD} $$
Taking Galois cohomology and using Hilbert's Theorem 90 we find
$$ \begin{CD}   
     \Qt @>[2]>> \Qt @>>> \HH^1(G,\Z/2\Z) @>>> 1.
   \end{CD} $$
Thus $\Qt/\Qts \simeq \HH^1(G,\Z/2\Z)$, and (\ref{ES1}) gives 
rise to an exact sequence 
$$\begin{CD}  \cC(\Q) @>{[2]}>> \cC(\Q) @>>> \Qt/\Qts. \end{CD} $$ 
It remains to identify the last map.

To this end, recall the construction of $\HH^1$: given an
exact sequence of $G$-modules 
$$ \begin{CD}
   0 @>>> A @>>> B @>f>> C @>>> 0, \end{CD} $$
we get a homomorphism $C^G \lra \HH^1(G,A)$ as follows:
for $c \in C^G$, pick a $b \in B$ such that $f(b) = c$
and then define the cocycle $x$ by $x(\sigma) = \sigma(b) - b$;
the image of $c$ is then the equivalence class of $x$.  

This provides us with the isomorphism $\Qt/\Qts \simeq \HH^1(G,\Z/2\Z)$:
given a coset $a\Qts$, pick a preimage $\sqrt{a} \in \bQ$, and 
then define the cocycle $x: G \lra \Z/2\Z$ by 
$x(\sigma) = \sigma(\sqrt{\alpha}\,)/\sqrt{\alpha}.$   

Next we study the connecting homomorphism 
$\delta: \cC(\Q)/2\cC(\Q) \lra \HH^1(\Z/2\Z)$. Let $P = (r,s) \in \cC(\Q)$. 
The points $Q = (x,y) \in \cC(\bQ)$ such that $2Q = P$ 
are given by
$$ Q = \begin{cases}
       (\sqrt{r+2)},s/\sqrt{r+2)}\,) & \text{if}\ r \ne -2, \\
       (0,2/\sqrt{-\Delta}\,)        & \text{if}\ r = -2. 
   \end{cases} $$
Via the homomorphism $\HH^1(\Z/2\Z) \simeq \Qt/\Qts$, the cocycle 
corresponding to $Q$ is identified with the coset
$$ \delta(P) = \begin{cases} 
         (r+2)\Qts & \text{if}\ r \ne -2, \\
      -\Delta \Qts & \text{if}\ r = -2. \end{cases} $$
Thus $\delta$ can be identified with the Weil map 
$\alpha: \cC(\Q) \lra \Qt/\Qts$, and in particular $\alpha$
is a group homomorphism with kernel $2 \cC(\Q)$.

\section{Heights}

For proving that $\cC(\Z)$ is finitely generated, we need more than
just the fact that $\cC(\Z)/2\cC(\Z)$ is finite. This missing piece
of information will be provided by the theory of heights.

\subsection{The Naive Height}

For rational numbers $x = \frac{r}{s}$ in lowest terms, 
we define 
$$H(x) = \max \{ |r|, |s|\};$$ 
note that $H(0) = 1$ and $H(x) \ge 1$ for all $x \in \Q$. 
The following lemma is easy to prove:

\begin{lem}\label{LeH}
For $x, y \in \Q$ we have 
\begin{enumerate}
\item $H(xy) \le H(x)H(y)$;
\item $H(x^2) = H(x)^2$;
\item $\frac1{2H(y)} Hx) \le H(x + y) \le 2H(x)H(y)$;
\item for any $c > 0$, the set of all $x \in \Q$ with height 
      $H(x) < c$ is finite. 
\end{enumerate}
\end{lem}

The lower bound in (3) follows from the upper bound upon replacing
$x$ by $x+y$ and $y$ by $-y$. 

Our next goal is the definition of the `naive height' $H(P)$ of 
rational points $P$ on Pell conics. For rational points 
$P = (x,y) \in \cC(\Q)$ on a conic $\cC: X^2 - \Delta Y^2 = 4$ 
put $H(P) = H(x)$. We clearly have

\begin{prop}
Let $\cC: X^2 - \Delta Y^2 = 4$ be a Pell conic. For a given 
constant $c > 0$, the set of all rational points $P \in \cC(\Q)$
with height $H(P) < c$ is finite.
\end{prop}

These rational points have a special form:

\begin{lem}
Let $(x,y) \in \cC(\Q)$ be a rational point on the Pell conic
$\cC: X^2 - \Delta Y^2 = 4$. Then there exist integers
$r,s,n$ such that $x = \frac{r}{n}$, $y = \frac{s}{n}$ and 
$\gcd(r,n) = \gcd(s,n) = 1$.
\end{lem}

\begin{proof}
Write $x = \frac{r}{n}$, $y = \frac{s}{m}$ with $r,s \in \Z$,
$m, n \in \N$ and $\gcd(r,n) = \gcd(s,m) = 1$. Then 
$r^2m^2 - \Delta s^2n^2 = 4m^2n^2$ shows that $n^2 \mid r^2m^2$,
and since $\gcd(r,n) = 1$, we find $n^2 \mid m^2$ and $n \mid m$. 

Thus $m = kn$ for some integer $k$. This gives
$r^2k^2 - \Delta s^2 = 4k^2n^2$, hence $k^2 \mid \Delta s^2$; 
since $k \mid m$ and $\gcd(s,m) = 1$ we conclude that $k^2 \mid \Delta$,
which implies that $k = 1$ if $\Delta \equiv 1 \bmod 4$ and $k \mid 2$
if $\Delta \equiv 0 \bmod 4$. In the latter case, $4k^2 \mid \Delta$
implies $4k^2 \mid r^2k^2$, hence $2 \mid r$; but this implies 
$k^2 \mid d$ and thus $k = 1$ as claimed.
\end{proof}

We also need some information on the height of the $Y$-coordinates.

\begin{lem}
Let $(x,y) \in \cC(\Q)$ with $y = \frac{s}{n}$; then 
$|\Delta| s^2 \le 4H(P)^2$.
\end{lem}

\begin{proof}
We have $|\Delta|s^2 \le \max \{r^2, 4n^2\} \le 4H(P)^2$.
\end{proof}

Now we claim   

\begin{prop}\label{PHt}
Let $Q \in \cC(\Q)$ be fixed. Then for all $P \in \cC(\Q)$ we have 
\begin{enumerate} 
\item $\frac14H(P)^2 \le H(2P) \le 4H(P)^2$; \smallskip
\item $\frac1c H(P) \le H(P + Q) \le cH(P)$ for $c = 5H(Q)$. 
\end{enumerate}
\end{prop}

\begin{proof}
For $P = (x,y)$ we have $2P = (x^2 - 2, xy)$, hence $H(2P) = H(x^2 - 2)$.
Lemma \ref{LeH}.(3) applied with $y = 2$ now proves the first claim.
For the proof of the second claim let $P = (x,y)$, $Q = (z,w)$ 
with $x = \frac{r}m$, $y = \frac{s}m$, $z = \frac{t}n$, $w = \frac{u}n$, 
and $\gcd(r,m) = \gcd(t,n) = 1$. Then
$P+Q = (\frac{xz+yw\Delta}2, \frac{xw+yz}2) = 
      (\frac{rt+su\Delta}{2mn}, \frac{ru+st}{2mn})$.

Clearly $2|mn| \le 2H(P)H(Q)$; thus it is sufficient to bound the
numerator. Here we find
\begin{align*}
   H(P+Q) & \le |r| \cdot |t|+|s|\sqrt{\Delta} \cdot |u| \sqrt{\Delta} \\
          & \le H(P)H(Q) + 4H(P)H(Q) = 5H(P)H(Q).
\end{align*}
Replacing $Q$ by $-Q$ shows that $H(P-Q) \le 5H(P)H(Q)$.  
Applying this result to $P+Q$ instead of $P$ shows that
$H(P) \le 5H(P+Q)H(Q)$, and this finally shows that 
$H(P+Q) \ge \frac15 H(P)$.
\end{proof}

\subsection{The Canonical Height}

The (naive) logarithmic height of a rational point $P \in \cC(\Q)$ 
is defined by $h_0(P) = \log H(P)$. Recall that 
\begin{itemize}
\item  $|h_0(2P) - 2h_0(P)| < \log 4$ for all $P \in \cC(\Q)$;
\item  given $Q \in \cC(\Q)$, put $c  = h_0(Q) + \log 5$; then 
       $h_0(P+Q) \le h_0(P) + c$ for every $P \in \cC(\Q)$.
\end{itemize}

Now let us define a function $h: \cC(\Q) \lra \R_{\ge 0}$
by putting
$$ h(P) = \lim_{n \to \infty} \frac{h_0(2^nP)}{2^n}. $$

In order to see that this definition makes sense we have to 
check that the sequence $\{2^{-n} h_0(2^nP)\}$ is Cauchy.

We know that $|h_0(2Q) - 2h_0(Q)| \le \log 4$; then $n > m \ge 0$ 
implies 
\begin{align*}
|2^{-n}h_0(2^nP) - 2^{-m}h_0(2^mP)| & = 
   \Big|\sum_{j=m}^{n-1} (2^{-j-1}h_0(2^{j+1}P) - 2^{-j}h_0(2^jP)) \Big| \\
   & \le \sum_{j=m}^{n-1} 2^{-j-1}|h_0(2^{j+1}P) - 2h_0(2^jP)| \\
   & \le \sum_{j=m}^{n-1} 2^{-j-1} \log 4 < 2^{-m} \log 4.
\end{align*}
Since this expression can be made arbitrarily small by choosing
$m$ sufficiently large, the sequence is Cauchy, and $h(P)$ is defined.
Taking $m = 0$ in the inequality above and letting $n \lra \infty$ proves 

\begin{prop}
For all $P \in \cC(\Q)$, we have $|h(P) - h_0(P)| \le \log 4$.
\end{prop}

This immediately implies 

\begin{prop}
Let $\cC: X^2 - \Delta Y^2 = 4$ be a Pell conic. For a given 
constant $c > 0$, the set of all rational points $P \in \cC(\Q)$
with canonical height $h(P) < c$ is finite.
\end{prop}
   
Now we can easily derive the basic properties of the canonical height:

\begin{thm}
The canonical height $h: \cC(\Q) \lra \R_{\ge 0}$ on the Pell conic
$\cC: X^2 - \Delta Y^2 = 4$ has the following
properties:
\begin{enumerate}
\item $h(T) = 0$ if and only if $T \in \cC(\Q)_\tors$;
\item $h(2P) = 2h(P)$;
\item $h(P+Q) \le h(P) + h(Q)$;
\item $h(P) + h(Q) \le h(P-Q) + h(P+Q) \le 2h(P) + 2h(Q)$;
\item the square of the canonical height satisfies the parallelogram
      equality 
      $$h(P-Q)^2 + h(P+Q)^2 = 2h(P)^2 + 2h(Q)^2$$
\end{enumerate}
for all $P, Q \in \cC(\Q)$.
\end{thm}

\begin{proof}
\begin{enumerate}
\item If $T$ is a torsion point, then $h_0(T^k)$ attains only 
      finitely many values, hence is bounded; this implies that
      $h(T) = 0$.

      Now assume that $h(T) = 0$. Then $h(kT) = k\cdot h(T)$ for all
      $k \ge 1$. Since $|h(P) - h_0(P)|$ is bounded, the naive heights
      of the points $kT$ are bounded. But there are only finitely
      many points with bounded height, hence $\{kT: k \in \N\}$ is
      finite, and this implies that $T$ is a torsion point. 

\item Directly from the definition we get
      $$ h(2P) = \lim_{n \to \infty} \frac{h_0(2^{n+1}P)}{2^n}
               = 2 \lim_{n \to \infty} \frac{h_0(2^{n+1}P)}{2^{n+1}}
               = 2h(P). $$

\item  Now recall that $h_0(P+Q) \le h_0(P) + h_0(Q) + \log 2$; 
       this implies
       \begin{align*}
       h(P+Q) & = \lim_{n \to \infty} \frac{h_0(2^n(P+Q))}{2^n} \\
              & \le  \lim_{n \to \infty} \Big( \frac{h_0(2^nP)}{2^n}
                      + \frac{h_0(2^nQ)}{2^n} + \frac{\log 2}{2^n} \Big) \\
              & = h(P) + h(Q). 
       \end{align*}
\item Replacing $Q$ by $-Q$ shows that $h(P-Q) \le h(P) + h(Q)$, 
      and adding these inequalities yields
      $$ h(P+Q) + h(P-Q) \le 2h(P) + 2h(Q). $$
      Applying this inequality to $P-Q$ and $P+Q$ instead of $P$ and $Q$
      yields
      $$ h(P) + h(Q) \le h(P+Q) + h(P-Q), $$
      where we have used $h(2P) = 2h(P)$ and $h(2Q) = 2h(Q)$.
\item Let us return to $h(P+Q) \le h(P) + h(Q)$; replacing 
      $P$ by $P-Q$ yields $h(P-Q) \ge h(P) - h(Q)$. Similarly, 
      $h(P+Q) \ge h(P) - h(Q)$. Squaring and adding yields 
      $h(P+Q)^2 + h(P-Q)^2 \ge 2h(P)^2 + 2h(Q)^2$. 

      Replacing $P$ and $Q$ by $P+Q$ and $P-Q$ shows
      $4h(P)^2 + 4h(Q)^2 = h(2P)^2 + h(2Q)^2 \ge 2h(P+Q)^2 + 2h(P-Q)^2$,
      that is, $h(P+Q)^2 + h(P-Q)^2 \le 2h(P)^2 + 2h(Q)^2$. 

      These two inequalities imply the desired equality.
\end{enumerate}
This concludes the proof.
\end{proof}

As a corollary we note:

\begin{cor}
We have $h(mP) = m h(P)$ for all $m \ge 1$.
\end{cor}

\begin{proof}
Put $P = mQ$ in the parallelogram equality.
\end{proof}

It is not hard to give explicit formulas for the canonical height 
of rational points on Pell conics:

\begin{prop}
The canonical height of $P = (x,y) \in \cC(\Q)$, where $\cC$ is the 
Pell conic given by $X^2 - \Delta Y^2 = 4$ with $\Delta > 0$, is  
$h(P) = \log \frac{|r|+|s|\sqrt{\Delta}}2$, where 
$x = \frac{r}{n}$, $y = \frac{s}{n}$ with $(r,n) = (s,n) = 1$.
\end{prop}

\begin{proof}
Observe that $2P = (\frac{r^2-2n^2}{n^2},\frac{rs}n)$ with
$(r^2-2n^2,n^2) = 1$, hence $H(2P) = r^2 - 2n^2$. Also note that
$r^2 - 2n^2 = n^2[(\frac{r+s\sqrt{\Delta}}{2n}\,)^2 + 
                  (\frac{r-s\sqrt{\Delta}}{2n}\,)^2].$
By induction, we conclude that for $k = 2^m$ and $r, s > 0$ we have
\begin{align*}
h(P) & = \lim_{k \to \infty} \frac{h_0(kP)}{k} = \lim_{k \to \infty}
         \frac1k  \log  n^k \Big[\Big(\frac{r+s\sqrt{\Delta}}{2n}\,\Big)^k +
                         \Big(\frac{r-s\sqrt{\Delta}}{2n}\,\Big)^k\Big] \\
     & = \log n + \lim_{k \to \infty}
         \frac1k  \log  \Big(\frac{r+s\sqrt{\Delta}}{2n}\,\Big)^k
       = \frac{r + s \sqrt{\Delta}}2,
\end{align*}
where we have used that $-1 < \frac{r-s\sqrt{\Delta}}{n} < 1$.
The other cases (e.g. $r > 0$, $s < 0$) are handled similarly.
\end{proof}

There is an even simpler formula if $\Delta < 0$:

\begin{prop}
The canonical height of $P = (x,y) \in \cC(\Q)$, where $\cC$ is the 
Pell conic given by $X^2 - \Delta Y^2 = 4$ with $\Delta < 0$, is  
$h(P) = \log n$, where $x = \frac{r}{n}$, $y = \frac{s}{n}$ with 
$n > 0$ and $(r,n) = (s,n) = 1$.
\end{prop}

\begin{proof}
We have $2^j P = (x_j, y_j)$, where $(x_k)$ is the sequence 
defined recursively by $x_1 = x$ and $x_{j+1} = x_j^2 - 2$. 
Clearly we have $|x_j| < 2$ for all $j \ge 1$, so the sequence 
is bounded. 

Assume that $|x_k| > 1$ for some $k$; we claim that there is a $j > 0$
such that $|x_{k+j}| < 1$. If not, we may assume that 
$x_k > 1$ (the case $x_k < -1$ is treated in an analogous way);
then $|x_{k+1}| > 1$ implies $x_{k+1} > 1$. On the other hand,
it is easily seen that in this case $x_{k+1} < x_k$. Thus
$1 < x_{k+j} > 1$ for all $j \ge 0$, hence the sequence converges, 
and we have $1 \le \lim x_j \le x_k < 2$; but the only possible limits
are the roots of the equation $0 = x^2 - x - 2 = (x+1)(x-2)$,
that is, $x = -1$ or $x = 2$: contradiction.
 
Thus there are infinitely many $x_k$ with $|x_k| < 1$; 
if we write $x = \frac{r}{n}$ with $n > 0$, then $x_k = r'/n^k$,
hence $H(2^kP) = H(x_k) = n^{2^k}$. We know that $\log 2^{-j} H(2^j P)$
converges to $h(P)$, hence so does the subsequence
$ 2^{-k} \log H(x_k) = \log n$.
\end{proof}

Finally, let us look at the heights of points on descendants.
If $P = (r,s)$ is a rational point on the descendant 
$\cT_a: ar^2 - bs^2 = 4$ with $ab = \Delta > 0$ and $a > 0$, then 
$Q = (ar^2-2, rs) \in \cC(\Q)$, and now Lemma \ref{LeH} implies
$$\frac1{4a} H(r)^2 \le \frac14 H(ar^2) \le H(Q) \le 4H(ar)^2 \le 4aH(r^2). $$

We have proved 

\begin{prop}
If If $P = (r,s)$ is a rational point on the descendant 
$\cT_a: ar^2 - bs^2 = 4$ with $ab = \Delta > 0$ and $a > 0$,
then $Q = (ar^2-2, rs) \in \cC(\Q)$ satisfies
$$ \frac1{4a} H(P)^2 \le H(Q) \le 4aH(P)^2. $$
\end{prop}

\section{The Theorem of Mordell-Weil}

The Theorem of Mordell-Weil states that the group of rational 
points on an elliptic curve defined over $\Q$ is finitely
generated. Its analog for conics says that the group of
integral points on a Pell conic is finitely generated
(more generally it can be shown that the group of $S$-integral
points on a Pell conic is finitely generated if $S$ is finite).

\subsection{Mordell-Weil}

We now show that $\cC(\Z)$ is finitely generated. The following
result is the abstract kernel of the proof:

\begin{thm}
Let $G$ be an abelian group such that $G/2G$ is finite. Assume that
there exists a function $h: G \lra \R_{\ge 0}$ with the following
properties:
\begin{enumerate}
\item For every $c > 0$, the set $\{g \in G: h(g) < c\}$ is finite;
\item We have $h(2g) = 2h(g)$ for all $g \in G$;
\item $h(g-g')^2 + h(g+g')^2 = h(g)^2 + h(g')^2$ for all $g, g' \in G$.
\end{enumerate}
Then $G$ is finitely generated.
\end{thm}

\begin{proof}
Let $\Gamma$ be a set of representatives of the finitely 
many cosets of $G/2G$. Then each $g \in G$ can be written 
as $g - \gamma = 2g'$ for some $\gamma \in \Gamma$ and a 
$g' \in G$. Put $c = \max \{h(\gamma): \gamma \in \Gamma\}$.

Now let $\Omega$ denote the subgroup of $G$ generqated by 
all the elements of $\Gamma$ and the (finitely many) elements 
$g \in G$ with $h(g) \le c$. We claim that $G = \Omega$.

If not, then let $g$ be an element in $G$ with minimal height
such that $g \notin \Omega$; observe that $h(g) > c$.
We can write $g - \gamma = 2g'$ for some $\gamma \in \Gamma$
and $g' \in G$, and find
$$ 4h(g')^2 = h(g-\gamma)^2 = 2h(g)^2 + 2h(\gamma)^2 - h(g+\gamma)^2 
            \le 2h(g)^2 + 2c^2 < 4h(g)^2.$$
Thus $h(g') < h(g)$, hence $g' \in \Omega$. But then so is
$g = 2g' + \gamma$: contradiction.
\end{proof}

Applying this to our situation we find

\begin{cor}
Let $\cC: X^2 - \Delta Y^2 = 4$ be a Pell conic. Then 
the group $\cC(\Z)$ is finitely generated, that is, 
$\cC(\Z) \simeq \cC(\Z)_\tors \oplus \Z^r$ for some finite
group $\cC(\Z)_\tors$ and some integer $r \ge 0$ called the rank of
$\cC$. Moreover, $\im \alpha = 2^{r+1}$.
\end{cor}

The torsion group of $\cC(\Q)$ is easy to determine: torsion 
points $(x,y)$ have integral coordinates, and we have $y = 0$ 
or $y = \pm 1$. In fact, if $k \ge 2$ is an integer and 
$P \ne N = (2,0)$ a rational point on $\cC$ with $kP = N$, then 
$\Q(\zeta_k) \subseteq \Q(\sqrt{\Delta}\,)$. Thus 
$$\cC(\Q)_\tors = \begin{cases}
        \{(\pm 2, 0), (\pm 1, \pm 1)\} & \text{if} \ \Delta = -3, \\
        \{(\pm 2, 0), (0, \pm 2)\}     & \text{if} \ \Delta = -4, \\
        \{( \pm 2,0)\}                 & \text{otherwise} \end{cases} $$

\section{Selmer and Tate-Shafarevich Groups}

The subset of curves $\cT(a): ar^2 - bs^2 = 4$ with a rational
point corresponds to a subgroup $\Sel_2(\cC)$ of $\Qt/\Qts$ 
called the $2$-Selmer group of $\cC$; we have already shown 
that if $\cT_a$ and $\cT_{a'}$ have a rational point, then
so does $\cT{a''}$, where $aa' = a''k^2$ for some positive 
and squarefree integer $a'' \mid \Delta$. The same argument
shows that the curves $\cT_a$ with an integral point form
a group $W_2(\cC)$, which is clearly a subgroup of $\Sel_2(\cC)$
isomorphic to $\im \alpha$. The $2$-part of the Tate-Shafarevich 
group $\TS_2(\cC)$ is then defined by the exact sequence
\begin{equation}\label{ESTS}
    \begin{CD} 
   1 @>>> W_2(\cC) @>>> \Sel_2(\cC) @>>> \TS_2(\cC) @>>> 1. 
   \end{CD} 
\end{equation}
In this section, we shall study these groups.

\subsection{The $2$-Selmer Group}

\begin{prop}\label{PLE}
The first descendant $\cT(a): ax^2 - by^2 = 4$, where $ab = \Delta$ 
and $a > 0$, has a rational point if and only if $(a/q) = (-b/p) = +1$ 
for all odd primes $p \mid a$ and $q \mid b$.
\end{prop}

\begin{proof}
Legendre's theorem states that the 
ternary quadratic form $ax^2 + by^2 + cz^2$, where $a, b, c \in \Z$ 
are coprime and squarefree, represents $0$ over the integers if and 
only if it represents $0$ over the reals and over the fields $\Z/p\Z$, 
where $p$ runs through the odd primes dividing $abc$.
\end{proof}

Note that $ax^2 - by^2 = 4$ has rational solutions if and only if
$X^2 = aY^2 + \Delta Z^2$ has integral solutions, hence the
criteria in Proposition \ref{PLE} are equivalent to 
$(\frac{a,\Delta}{p}) = +1$ for all odd primes $p \mid \Delta$;
since $\Delta > 0$, the Hilbert symbol at $\infty$ is trivial;
finally, $(\frac{a,\Delta}{p}) = +1$ for all odd primes 
$p \nmid \Delta$, and now the product formula implies that we 
have $(\frac{a,\Delta}{2}) = +1$ as well. This shows

\begin{cor}
The first descendant $\cT(a): ax^2 - by^2 = 4$, where 
$ab = \Delta$ and $a > 0$, has a rational point if and 
only if $(\frac{a,\Delta}{p}) = +1$ for all primes $p$. 
\end{cor}

\medskip\noindent {\bf Example.}
Consider $\cC(\Z)$ for $\cC: x^2 - 1045y^2 = 4$. The associated 
descendants with an integral point $(r,s)$ and the corresponding 
point $(x,y)$ on $\cC$ are given below:
$$ \begin{array}{r|r|c|c}
   a & \cT_a(\cC) \quad  & (r,s) & (x,y) \\ \hline
   1 &     r^2 - 1045s^2 = 4 & (2,0) & (2,0) \\
   5 &    5r^2 -  209s^2 = 4 & (\frac73,\frac13)  &   \\
  11 &   11r^2 -   95s^2 = 4 & (3,1) &  (97,3) \\
  19 &   19r^2 -   55s^2 = 4 &  --   &   \\
  55 &   55r^2 -   19s^2 = 4 & (\frac47, \frac67) & \\
  95 &   95r^2 -   11s^2 = 4 &  --   &   \\
 209 &  209r^2 -    5s^2 = 4 &  --   &   \\
1045 & 1045r^2 -     s^2 = 4 &  --   &   \end{array} $$
Thus $\Sel_2(\cC)$, viewed as a subgroup of $\Qt/\Qts$,
is isomorphic to $\langle 5, 11 \rangle$; moreover
$W_2(\cC) = \langle 11 \rangle$, and the nontrivial
element of $\TS(\cC)[2] \simeq \Z/2\Z$ is generated by
$\cT_5$.

\subsection{R\'edei}
Recall from \cite{pell1} that a factorization of the discriminant 
$\Delta = \disc k$ into discriminants $\Delta = \Delta_1 \Delta_2$ 
is called a splitting of the second kind if 
$(\Delta_1/p_2) = (\Delta_1/p_2) = +1$ for all primes 
$p_i \mid \Delta_i$.

\begin{prop}
Assume that $\Delta$ is a product of positive prime discriminants. 
Then the factorization $\Delta = \Delta_1 \Delta_2$ is a splitting 
of the second kind if and only if the descendant 
$\Delta_1 X^2 - \Delta_2Y^2 = 4$ is everywhere locally solvable.
\end{prop}  

Thus R\'edei's group structure on splittings of the second kind
induces a group structure on $\Sel_2(\cC)$ that coincides with ours.
Since there are exactly $e_4+1$ independent splittings of the second
kind (including the trivial factorization $\Delta = 1 \cdot \Delta$),
this shows that R\'edei's results imply that 
$\# \Sel_2(\cC) = 2 \# \Cl^+(k)^2/\Cl^+(k)^4$. 
 
In general, however, the $C_4$-decompositions and the 
first descendants in the Selmer group are not related.
Consider e.g. the example $d = 12369 = 3 \cdot 7 \cdot 19 \cdot 31$;
here the elements in the Selmer group and the corresponding rational
points are given by 
\begin{align*}
     r^2 - 12369s^2 & = 4 &  (r,s) & = (2,0) \\
    7r^2 -  1767s^2 & = 4 &  (r,s) & \ts = (\frac{32}5, \frac25) \\ 
  589r^2 -    21s^2 & = 4 &  (r,s) & \ts = (\frac14, \frac54) \\
 4123r^2 -     3s^2 & = 4 &  (r,s) & \ts = (\frac1{32}, \frac3{32})
\end{align*}
The $C_4$-decompositions, on the other hand, are
$\Delta = 1 \cdot 12369$ and  $\Delta = 93 \cdot 133$,
and the equation $93 x^2 + 133y^2 = z^2$ has the solution
$(x,y,z) = (6,1,59)$.

\subsection{The $2$-Part of the Tate-Shafarevich Group}

Consider the descendant $\cT(a): ar^2 - bs^2 = 4$; we know that
$\cT(a) \in \Sel_2(\cC)$ if and only if $(\frac{a,d}{p}) = +1$
for all places $p$ (observe that $(\frac{a,d}{p}) = +1$ for all 
primes $p \nmid d$).

Consider the map $\cl:\Sel_2(\cC) \lra \Cl^+(k)[2]$ sending
$\cT(a)$ to the ideal class generated by the ambiguous ideal 
$\fa$ with norm $a$; clearly $\ker \cl = W(\cC)$. By Hilbert's 
genus theory (see \cite{LRL}), we know that ideal classes coming 
from the Selmer group are squares, so the map above is actually 
a homomorphism $\Sel_2(\cC) \lra \Cl^+(k)^2 \cap \Cl^+(k)[2]
 = \Cl^+(k)^2[2]$. Conversely, an ideal class in $\Cl^+(k)^2[2]$
is generated by an ambiguous ideal $\fa$ with norm $a \mid \Delta$, 
and since its class is a square, its character system is trivial, 
so the descendant $\cT(a)$ is in the Selmer group. 

\begin{thm}\label{TTS}
We have an exact sequence
$$ \begin{CD} 
   0 @>>> W(\cC) @>>> \Sel_2(\cC) @>>> \Cl^+(k)^2[2] @>>> 0. \end{CD} $$
In particular, $\TS_2(\cC) \simeq \Cl^+(k)^2[2]$.
\end{thm}

Observe that, for finite abelian groups $G$, we have the exact sequence
$$ \begin{CD} 
   1 @>>> G^2 \cap G[2] @>>> G^2 @>[2]>> G^4 @>>> 1 \end{CD} $$
showing that $G^2 \cap G[2] \simeq G^2/G^4$ (non-canonically
via duality), hence $\Cl^+(k)^2[2] \simeq \Cl^+(k)^2/\Cl^+(k)^4$.
Since this group can be made arbitrarily large, we find

\begin{cor}
For Pell conics $\cC: X^2 - \Delta Y^2 = 4$, the Tate-Shafarevich 
group $\TS_2(\C)$ can have arbitrarily large $2$-rank as $\Delta$ 
varies.
\end{cor}

\subsection{For Whom the Pell Tolls}

Let us now derive some results about Pell equations that follow
from studying the $2$-Selmer group.

\subsection*{Selmer Groups}

The following result connects the structure of the Selmer group
to various invariants studied in \cite{pell1}:

\begin{prop}
Let $\Delta$ be a discriminant not divisble by any prime
$\equiv 3 \bmod 4$, let $\cC: X^2 - \Delta Y^2 = 4$ be the 
corresponding Pell conic, and let $\gamma(\Delta)$ be the 
associated nondirected graph (see \cite{pell1}). Then the 
following claims are equivalent:
\begin{enumerate}
\item $\gamma(\Delta)$ is odd;
\item $\Sel_2(\cC) \simeq \Z/2\Z$;
\item $\TS(\cC)[2] = 0$;
\item $4$-rank $\Cl_2^+(k) = 0$.
\end{enumerate}
\end{prop}

\begin{proof}
The equivalence of the statements (2)-(4) follow from the 
exact sequence (\ref{ESTS}) and Theorem \ref{TTS}.

The fact that $\gamma(\Delta)$ is odd if and only if
none of the equations $ar^2 - bs^2 = 4$ has $\F_p$-rational
points for all primes $p$ was proved in \cite{pell1}.
\end{proof}

\subsection*{Nontrivial Tate-Shafarevich Groups}
If $a$ is a quadratic residue modulo a prime $p$, then 
we write $(\frac{a}{p})_4 = +1$ or $-1$ according as $a$ 
is a fourth power modulo $p$ or not. If $p \equiv 1 \bmod 8$,
then we define $(\frac{p}{2})_4 = (-1)^{(p-1)/8}$. We 
extend these residue symbols multiplicatively to composite 
denominators.

\begin{thm}\label{TnC}
Let $\Delta = p_1 \cdots p_n$ be a product of primes
$p_i \equiv 1 \bmod 4$. If $ab = \Delta$ and $ar^2 - bs^2 = 4$ 
has an integral solution, then the following conditions are 
satisfied:
\begin{enumerate}
\item $(a/q) = 1$ for all primes $q \mid b$;
\item $(b/p) = 1$ for all primes $p \mid a$;
\item $(b/a)_4 = +1$.
\end{enumerate}
\end{thm}

\begin{proof}
The first two assertions are clear and follow from the existence 
of a rational point.

For any prime $p \mid a$, we have $(-b/p)_4 = (2s/p)$; 
since $(-1/p)_4 = (2/p)$, this implies $(b/p)_4 = (s/p)$,  
hence $(b/a)_4 = (s/a)$. Now write $s = 2^js'$ with $s'$ odd;
then $(b/a)_4 = (s/a) = (2/a)^j$. If $j = 0$ or $j = 2$, we 
are done. The case $j = 1$ is impossible: putting $r = 2r'$
we find $a{r'}^2 - b{s'}^2 = 1$, which leads to a contradiction
modulo $4$ since $b{s'}^2 \equiv 1 \bmod 4$. Finally, if
$j \ge 3$, then dividing $ar^2 - bs^2 = 4$ through by $4$
and reducing modulo $8$ shows that $(2/a) = +1$.  
\end{proof}

There are similar results for even $\Delta$ not divisible by
primes $\equiv 3 \bmod 4$.

Let us now apply this result to the Pell equation $X^2 - pqY^2 = 1$,
where $p \equiv q \equiv 1 \bmod 4$. The first descendants
$pr^2 - qs^2 = \pm 1$ are not solvable in integers if $(p/q) = -1$,
so in this case we conclude that $X^2 - pqY^2 = -1$ is solvable.
Assume that $(p/q) = +1$. Then Theorem \ref{TnC} provides us with 
necessary conditions for the descendant $\cT_a$ to be solvable:
\begin{table}[ht!]\label{T1}
$$ \begin{array}{r|r|c}
    a & \text{equation} & \text{condition} \\ \hline 
  \rsm  1 &  X^2 - pqY^2 = 1 & \text{none} \\ 
  \rsm  p &  pX^2 - qY^2 = 1 & (q/p)_4 = 1 \\ 
  \rsm  q &  qX^2 - pY^2 = 1 & (p/q)_4 = 1 \\
  \rsm pq &  pqX^2 - Y^2 = 1 &    ? 
   \end{array} $$
\caption{Solvability Criteria for $\cT$}
\end{table}
Thus if $(p/q)_4 = (q/p)_4 = -1$, the negative Pell equation 
$X^2 - pqY^2 = -1$ must be solvable. If, say, $(p/q)_4 = -(q/p)_4$,
however, we do not get a precise result because Theorem \ref{TnC}
does not give us any condition for the solvability of $\cT_{pq}$.
For this, we have to dig deeper:

\begin{prop}\label{PDd}
If $\cT_{pq}: pxr^2 - s^2 = 1$ has an integral solution, where 
$p \equiv q \equiv 1 \bmod 4$ are primes with $(p/q) = 1$, then
$(p/q)_4 = (q/p)_4$.
\end{prop}

This implies the following result, parts of which were first 
proved by Scholz \cite{Sch34} using class field theory:

\begin{prop}\label{PSch}
Let $p \equiv q \equiv 1 \bmod 4$ be primes. If the conditions
$(*)$ are verified, the descendant $\cT_a$ of the Pell conic
$X^2 - pqY^2 = 1$ is solvable:
$$ \begin{array}{c|c}
   (*)  & a \\ \hline
   \rsm  (p/q) = -1  & pq \\
   \rsm  (p/q) = +1, (p/q)_4 = -1, (q/p)_4 = +1 & p \\
   \rsm  (p/q) = +1, (p/q)_4 = +1, (q/p)_4 = -1 & q \\
   \rsm  (p/q) = +1, (p/q)_4 = -1, (q/p)_4 = -1 & pq
   \end{array} $$
\end{prop}

Note that this implies e.g. that if $(p/q) = +1$, 
$(p/q)_4 = +1$ and $(q/p)_4 = -1$, then $\cT_q$ is an
element of the Selmer group without an integral point,
hence represents an element of order $2$ in $\TS(\cC)$.

The proof of Proposition \ref{PSch} presents no problems; thus
it remains to prove Proposition \ref{PDd}. This is done as
follows: factor the right hand side of $pqs^2 = r^2 + 1$ over
the Gaussian integers $\Z[i]$. Since $\gcd(r+i,r-i)$ divides
$2i$, and since $r$ is even, the factors $r+i$ and $r-i$
are coprime. Now observe that $p = \pi \bpi$ and $q = \rho \brho$
for $\pi, \rho \in \Z[i]$, where the bars denote the conjugates.  
Assume that $\pi$ and $\rho$ are primary, i.e., that 
$\pi \equiv \rho \equiv 1 \bmod 2+2i$. Then Unique Factorization
in $\Z[i]$ implies that $r+i = \eps \pi \rho \alpha^2$ for some
$\alpha \in \Z[i]$ and some unit $\eps \in \{\pm i, \pm 1\}$. 
Since $r+i \equiv i \bmod 2$, and since $\alpha^2 \equiv 1 \bmod 2$,
we have $\eps = \pm i$, and by subsuming the square $-1 = i^2$ into 
$\alpha$ if necessary we arrive at $r+i = i \pi \rho \alpha^2$.

If, from this equation, we subtract its conjugate and then
divide by $i$, we arrive at
$$ 2 = \pi \rho \alpha^2 - \bpi \brho \balpha^2. $$
Reducing modulo $\brho$ we find $[2/\brho] = [\pi/\brho][\rho/\brho]$,
where $[\,\cdot\,/\,\cdot\,]$ is the quadratic residue symbol in 
$\Z[i]$ (see \cite{LRL} for the necessary background). Then it is
known that $[2/\brho] = (2/q)$ and $[\rho/\brho] = (2/q)$, as well 
as $[\pi/\brho] = [\pi/\rho] = (p/q)_4 (q/p)_4$. This concludes
the proof of Proposition \ref{PDd}.

This allows us to complete Table \ref{T1}:

\begin{table}[ht!]\label{T2}
$$ \begin{array}{r|r|c}
    a & \text{equation} & \text{condition} \\ \hline 
  \rsm  1 &  X^2 - pqY^2 = 1 & \text{none} \\ 
  \rsm  p &  pX^2 - qY^2 = 1 & (q/p)_4 = 1 \\ 
  \rsm  q &  qX^2 - pY^2 = 1 & (p/q)_4 = 1 \\
  \rsm pq &  pqX^2 - Y^2 = 1 & (p/q)_4 (q/p)_4 = 1     
   \end{array} $$
\caption{Solvability Criteria for $\cT$}
\end{table}

In some sense, the solvability condition for the product
$\cT_{pq}$ of $\cT_p$ and $\cT_q$ is the `product' of the
conditions for $\cT_p$ and $\cT_q$; although we cannot make
this more precise at the moment, this observation often
helps to guess the right criteria.

Observe that the proof of Proposition \ref{PSch} is fully 
analogous to the calculations done in \cite{LFNC} for computing 
Tate-Shafarevich groups of elliptic curves connected to the
congruent number problem.


\begin{thebibliography}{99}

\bibitem[Bal1999]{Bal} N. Baldisserri, 
{\em The group of primitive quasi-Pythagorean triples} (Italian),
Rend. Circ. Mat. Palermo (2) {\bf 48} (1999), 299--308; cf. p. 
%

\bibitem[BS1996]{BS1} R.A. Beauregard, E.R. Suryanarayan,
{\em Pythagorean triples: the hyperbolic view},
College Math. J. 1996; cf. p. 
%

\bibitem[BS1997]{BS2} R.A. Beauregard, E.R. Suryanarayan,
{\em Arithmetic Triangles},
Math. Mag. {\bf 70} (1997), 105--115; cf. p. 
%

\bibitem[BS1999]{BS3} R.A. Beauregard, E.R. Suryanarayan,
{\em Integral Triangles},
Math. Mag. {\bf 72} (1999), 287--294; cf. p. 
%

\bibitem[Daw1994]{Daw} B. Dawson,
{\em The ring of Pythagorean triples},
Missouri J. Math. Sci. {\bf 6} (1994), 72--77; cf. p. 
%

\bibitem[Dic1920]{Dick} L.E. Dickson,
{\em History of the Theory of Numbers},
vol I (1920); vol II (1920); vol III (1923);
Chelsea reprint 1952; cf. p. 
%

\bibitem[Dic1930]{DiSN} L.E. Dickson,
{\em Studies in the Theory of numbers}, 
Chicago 1930; cf. p. 
%

\bibitem[Eck1984]{Eck} E. Eckert,
{\em The group of primitive Pythagorean triangles},
Math. Mag. {\bf 54} (1984), 22--27; cf. p. 
%

\bibitem[Gry1997]{Gry} A. Grytczuk, 
{\em Note on a Pythagorean ring},
Missouri J. Math. Sci. {\bf 9} (1997), 83--89; cf. p. 
%

\bibitem[Hla2000]{Hla} E. Hlawka,
{\em Pythagorean triples},
Number theory, Birkh\"auser, Basel (2000), 141--155; cf. p. 
%

\bibitem[Jue1896]{Juel} C. Juel,
{\em Ueber die Parameterbestimmung von Punkten auf
Curven zweiter und dritter Ordnung. Eine geometrische
Einleitung in die Theorie der logarithmischen und
elliptischen Funktionen}, Math. Ann. {\bf 47} (1896), 72--104; cf. p.
%

\bibitem[Lem2000]{LRL} F. Lemmermeyer,
{\em Reciprocity Laws. From Euler to Eisenstein},
Springer Verlag 2000; cf. p.
%

\bibitem[Lem2003a]{LFNC} F.~Lemmermeyer,
{\em Some families of non-congruent numbers},
Acta Arith. {\bf 110} (2003), 15--36
%

\bibitem[Lem2003b]{pell1} F. Lemmermeyer,
{\em Higher Descent on Pell Conics. I. From Legendre to Selmer},
preprint 2003; cf. p.
%

\bibitem[Lem2003c]{pell2} F. Lemmermeyer,
{\em Higher Descent on Pell Conics. II. Two Centuries of Missed Opportunities},
preprint 2003; cf. p.
%

\bibitem[Mar1962]{Mar} J. Mariani,
{\em The group of the pythagorean numbers},
Amer. Math. Mon. {\bf 69} (1962), 125--128; cf. p.
%

\bibitem[Mor1986]{Mor} J. Morita, 
{\em A transformation group of the Pythagorean numbers},
Tsukuba J. Math. {\bf 10} (1986), no. 1, 151--153; cf. p.
%

\bibitem[Nie1908]{Nie} B. Niewenglowski, 
{\em Note sur les equations $x^2 - ay^2 = 1$ et $x^2 - ay^2 = -1$},
Bull. Soc. Math. France {\bf 35} (1907), 126--131; cf. also
Wiadomi Mat. Warsaw {\bf 12} (1908), 1--26 (Polish); cf. p.
%

\bibitem[PS1997]{PrS} V. Prasolov, Y. Solovyev,
{\em Elliptic Functions and Elliptic Integrals},
Transl. Math. Monographs {\bf 170}, AMS 1997; cf. p. 
%

\bibitem[Sch1990]{Sch} N. Schappacher,
{\em D\'eveloppement de la loi de groupe sur une cubique},
S\'eminaire Th\'eor. Nombres, Paris 1988--1989, 159--184;
Progr. Math. {\bf 91} (1990); cf. p. 
%

\bibitem[Sch1839]{Schm} Th. Sch\"onemann,
{\em Ueber die Congruenz $x^2 + y^2 \equiv 1 \pmod p$},
J. Reine Angew. Math. {\bf 19} (1839), 93--112; cf. p. 
%

\bibitem[Sch1934]{Sch34} A.~Scholz,
{\em  \"Uber die L\"osbarkeit der Gleichung $t^2-Du^2 =-4$},
 Math. Z. {\bf 39} (1934),  95--111; cf. p. 
%

\bibitem[Sha2001]{Sha} P. Shastri, 
{\em Integral points on the unit circle},
J. Number Theory {\bf 91} (2001), 67--70; cf. p. 
%

\bibitem[Sta1896]{Stae} P. St\"ackel,
{\em Review JFM 27.0337.02}, 
Jahrbuch Fortschritte der Mathematik 
{\bf 27} (1896), p. 337; cf. p.  
%

\bibitem[Tan1996]{Tan} L. Tan,
{\em The group of rational points on the unit circle},
Math. Mag. {\bf 69} (1996), 163--171; cf. p. 
%

\bibitem[Tau1970]{Tau} O. Taussky,
{\em Sums of squares},
Amer. Math. Monthly {\bf 77} (1970), 805--830; cf. p.
%

\bibitem[Tur1915]{Tu15} E. Turri\`ere, 
{\em Le probl\`eme de {\it Jean de Palerme} et {\it de L\'eonard de Pise}},
Ens. Math. {\bf 17} (1915), 315--324; cf. p.
%

\bibitem[Tur1916]{Tu16} E. Turri\`ere, 
{\em  Notions d'arithmog\'eom\'etrie},
Ens. math. {\bf 18} (1916), 81--110, 397--428; cf. p.
%

\bibitem[Tur1917]{Tu17} E. Turri\`ere, 
{\em  Notions d'arithmog\'eom\'etrie},
Ens. math. {\bf 19} (1917), 159--191, 233--272; cf. p.
%

\bibitem[Tur1918]{Tu18} E. Turri\`ere, 
{\em  Notions d'arithmog\'eom\'etrie},
Ens. math. {\bf 20} (1918), 161--174; cf. p.
%

\bibitem[VY1910]{VY} O. Veblen, J.W. Young,
{\em Projective Geometry I},
Ginn \& Co. 1910; cf. p.
%

\bibitem[Woi2001]{Woi} M. Wojtowicz,
{\em Algebraic structures on some sets of Pythagorean triples. II},
Missouri J. Math. Sci. {\bf 13} (2001), 17--23; cf. p. 
%

\bibitem[ZZ1991]{ZZ} P. Zanardo, U. Zannier,
{\em The group of Pythagorean triples in number fields},
Ann. Mat. Pura Appl. (4) {\bf 159} (1991), 81--88; cf. p. 
%

\end{thebibliography}
\end{document}